\documentclass{conm-p-l}

\usepackage{amsmath}
\usepackage{amsfonts}
\usepackage{amssymb}
\usepackage{amsthm}
\usepackage[all]{xy}

\newtheorem{lem}{Lemma}[section]

\newtheorem{thm}[lem]{Theorem}
\newtheorem{conj}[lem]{Conjecture}
\newtheorem{Defn}[lem]{Definition}
\newtheorem{Ex}[lem]{Example}
\newtheorem{Question}[lem]{Question}
\newtheorem{Problem}[lem]{Problem}

\newenvironment{ex}{\begin{Ex}\rm}{\end{Ex}}
\newenvironment{question}{\begin{Question}\rm}{\end{Question}}

\newtheorem{llem}{Lemma}[subsection]
\newtheorem{tthm}[llem]{Theorem}
\newtheorem{Property}[llem]{}
\newenvironment{property}{\begin{Property}\rm}{\end{Property}}

\newcommand{\ideal}[1]{\mathfrak{#1}}

\newcommand{\m}{\ideal{m}}
\newcommand{\n}{\ideal{n}}
\newcommand{\p}{\ideal{p}}
\newcommand{\q}{\ideal{q}}
\newcommand{\s}{\ideal{s}}
\newcommand{\pp}{\ideal{P}}
\newcommand{\qqq}{\ideal{Q}}
\newcommand{\mm}{\ideal{M}}

\newcommand{\x}{\mathbf{x}}
\newcommand{\rank}{\mathrm{rank}}	
\newcommand{\Ht}{\mathrm{ht}\,}	

\newcommand{\len}{\mathrm{len}}
\newcommand{\tor}{\text{Tor}}

\begin{document}

\title[Symbolic Powers of Prime Ideals]{On Symbolic Powers of Prime Ideals}
\author{Sean Sather-Wagstaff}
\address{Department of Mathematics, University of Illinois, 273 Altgeld Hall,
1409 West Green Street, Urbana IL, 61801, USA}
\email{ssather@math.uiuc.edu}
\thanks{This research was conducted while the author was an NSF Mathematical 
Sciences Postdoctoral Research Fellow}
\copyrightinfo{2003}
{American Mathematical Society}  

\begin{abstract}
Let $(R,\m)$ be a regular local ring with prime ideals $\p$ and $\q$ such 
that $\sqrt{\p+\q}=\m$ and $\dim(R/\p)+\dim(R/\q)=\dim(R)$.  It has been 
conjectured by Kurano and Roberts that $\p^{(m)}\cap\q\subseteq\m^{m+1}$ 
for all positive integers $m$.  We discuss this conjecture and related 
conjectures.  In particular, we prove that this conjecture holds for all 
regular local rings if and only if it holds for all localizations of 
polynomial algebras over complete discrete valuation rings.
In addition, we give an example showing that a certain 
generalization to nonregular rings does not hold.
\end{abstract}

\maketitle






\section{Introduction} \label{sec:intro}
All rings in this paper are 
assumed commutative and Noetherian, and all modules unital.

Recent years have seen a renewed interest in the properties of symbolic 
powers of prime ideals, especially in regular local rings.  
For instance, see~\cite{eisenbud:essfi, hochster:csopi, huckaba:pihsad, 
huneke:ssrlr, kurano:pimsppi, sather:mdium, sather:dicmr, sather:isppi}.
Often these 
properties, or conjectured properties, are related to the behavior of the 
symbolic powers $\p^{(m)}$ of a prime ideal $\p$ with respect to the 
the maximal ideal $\m$.  These properties frequently
emerge from seemingly unrelated questions in algebra and geometry.  
For instance, we have the following.  (``SP'' stands for ``Symbolic 
Powers''.)  


\begin{conj} \label{conj:sp}
Let $(R,\m)$ be a regular local 
ring with prime ideals $\p$ and $\q$ such that $\sqrt{\p+\q}=\m$ and 
$\dim(R/\p)+\dim(R/\q)=\dim(R)$.
\begin{description}
\item[\rm(SP-1)] $\p^{(m)}\cap\q\subseteq\m^{m+1}$ for all $m\geq 1$.
\item[\rm(SP-2)]  $\p^{(m)}\cap\q^{(n)}\subseteq\m^{m+n}$ for all $m,n\geq 1$.
\end{description}
\end{conj}

(SP-1) is a conjecture of Kurano-Roberts~\cite[3.1]{kurano:pimsppi}, 
motivated by their work on Serre's Positivity Conjecture;  
we briefly discuss the connection between these two conjectures 
in Section~\ref{sec:back} below.
Kurano and Roberts verify (SP-1) in~\cite[3.2]{kurano:pimsppi}
for rings containing a field, and the 
current author gives a different proof in~\cite{sather:dicmr}.  
(SP-2) was conjectured and verified for rings containing a field
by the current author in~\cite[1.6, 2.8]{sather:isppi}.  Each 
conjecture is open in mixed characteristic.  In our main theorem
we show that the ramified case 
of each conjecture follows from the unramified case.  

The current author (op.~cit.) has also made the following conjecture, 
which is fundamentally related to Conjecture~\ref{conj:sp}.  
First, we explain a notational convention.  
In general, regular rings are denoted by $R$, and their 
prime ideals by Fraktur 
characters, e.g., $\p$ and $\q$.  
Rings that need not be regular are denoted by $A$
and their prime ideals 
by Roman characters, e.g., $P$ and $Q$; the main
exception being the maximal ideal $\n$.
We give the definitions of technical terms in the conjecture after the 
statement.
(``ID'' stands for ``intersection dimension''.)

\begin{conj} \label{conj:id}
Let $(A,\n)$
be a quasi-unmixed local ring with prime ideals $P$ and 
$Q$ such that $\sqrt{P+Q}=\n$.  
\begin{description}
\item[\rm(ID-1)]  If $e(A)=e(A_P)$ and $A/P$ is analytically unramified, then 
\[ \dim(A/P)+\dim(A/Q)\leq\dim(A).\]
\item[\rm(ID-2)]  If $e(A)<e(A_{P})+e(A_{Q})$ and both $A/P$ and $A/Q$ are 
analytically unramified,  then 
\[ \dim(A/P)+\dim(A/Q)\leq\dim(A).\]
\end{description}
\end{conj}

In the statement, $e(A)$ is the Hilbert-Samuel multiplicity of the local 
ring $A$ with respect to its maximal ideal.  
A ring $A$ is \textit{quasi-unmixed} if the 
completion $A^*$ is equidimensional, i.e., for every minimal prime ideal 
$P^*$ of $A^*$ we have $\dim(A^*/P^*)=\dim(A)$.  A theorem of 
Ratliff~\cite{ratliff:qldafccpiII} says that $A$ is quasi-unmixed if and 
only if it is equidimensional and universally catenary.  
In particular, when $A$ is excellent, then $A$ is quasi-unmixed if and 
only if it is equidimensional.
The quotient ring $A/P$ is \textit{analytically 
unramified} if $(A/P)^*=A^*/PA^*$ is reduced, i.e., the 
ideal $PA^*$ is an intersection of prime ideals of $A^*$.  This is 
automatic if $A$ is excellent.

For rings that contain a field, the current author has verified both
(ID-1)~\cite[3.2]{sather:mdium} and (ID-2)~\cite[2.2]{sather:isppi}.
As with Conjecture~\ref{conj:sp}, each 
conjecture (ID-$i$) is open in mixed characteristic.  
In our main theorem below, we show that the general case of each
conjecture follows from the case of a hypersurface over an unramified
regular local ring.  


Some relations between our conjectures are straightforward to verify.
We summarize them in the following diagram 
\begin{equation} \label{diagram}
\begin{gathered}
\xymatrix{
\text{(ID-2)} \ar@{=>}[r] \ar@{=>}[d] & \text{(SP-2)} \ar@{=>}[d] \\
\text{(ID-1)} \ar@{=>}[r] & \text{(SP-1)}  } 
\end{gathered}
\end{equation}
and defer explanation to Section~\ref{sec:back}.
In our main theorem, which we now state,
 we show that the converses of the horizontal
implications in this diagram hold.  

\medskip

\noindent
\textbf{Main Theorem.}
\textit{Fix $i=1$ or $2$.  With the notations of Conjectures~\ref{conj:sp}
and~\ref{conj:id}, the following conditions are equivalent.
\begin{itemize}
\item[(a)]  \emph{(SP-$i$)} holds for all regular local rings.
\item[(b)]  \emph{(SP-$i$)} holds for all unramified regular local rings
that are essentially of finite type over a discrete valuation ring of mixed
characteristic.
\item[(c)]  \emph{(ID-$i$)} holds for all quasi-unmixed local rings.
\item[(d)]  \emph{(ID-$i$)} holds for all hypersurfaces over unramified 
regular local rings
that are essentially of finite type over a discrete valuation ring of mixed
characteristic.
\end{itemize}}

Here we summarize the contents of this paper.  In Section~\ref{sec:back}
we give some background material.  In Subsection~\ref{subsec:history} 
we discuss
the history of our conjectures.  In Subsection~\ref{subsec:imply} we 
verify some implications between the conjectures.  In 
Subsection~\ref{subsec:tools} we include definitions and tools relevant
to our results.  In Subsection~\ref{sec:ta} we discuss
the technical assumptions in the conjectures.  In 
Section~\ref{sec:main}, we prove our main theorem
and verify a weak version of (ID-2) in mixed characteristic.  
Finally, in Section~\ref{sec:ex} we give  
an example to show that a potential 
generalization of (SP-1)
does not hold.

\section{Background} \label{sec:back}

We begin this section by discussing the connection between Serre's 
Positivity Conjecture and (SP-1), and the other evidence supporting 
our conjectures.

\subsection{Evidence for the conjectures} \label{subsec:history}

Let $(R,\m)$ be a regular 
local ring of dimension $d$, and $\p,\q$ prime ideals of $R$ such that 
$\sqrt{\p+\q}=\m$.  
Serre~\cite[Th\'{e}or\`{e}me 3 of Chapitre V]{serre:alm} proved that
\[ \dim(R/\p)+\dim(R/\q)\leq d.\]
We shall refer to this result as Serre's Intersection Theorem.
Serre defined the 
\textit{intersection multiplicity} of $R/\p$ and $R/\q$ by the formula
\[ \chi(R/\p,R/\q)=\sum_{i=0}^d(-1)^i\len(\tor_i(R/\p,R/\q)) \]
where $\len(T)$ is the length of the $R$-module $T$.  
For $R$ unramified, Serre 
proved the following:  (i) $\chi(R/\p,R/\q)\geq 0$, 
and (ii) $\chi(R/\p,R/\q)=0$ if and only if $\dim(R/\p)+\dim(R/\q)<d$.  He 
conjectured that these results hold even when $R$ is ramified.  
Gillet-Soul\'{e}~\cite[Th\'{e}or\`{e}me 1]{gillet:knmi} and 
Roberts~\cite[Theorem 1]{roberts:vimpc} independently proved the 
Vanishing Conjecture:  if 
$\dim(R/\p)+\dim(R/\q)<d$, then $\chi(R/\p,R/\q)=0$.  
Gabber~\cite{berthelot:ava}, \cite{hochster:nimrrlr}, 
\cite{roberts:rdsmcgpnc} proved the 
Nonnegativity Conjecture:  $\chi(R/\p,R/\q)\geq 0$.  The Positivity 
Conjecture is the converse of the Vanishing Conjecture, and is still open.

The following theorem of Kurano-Roberts~\cite[3.2]{kurano:pimsppi}
shows the connection between 
the Positivity Conjecture and (SP-1).  
Essentially, 
it says that, for rings that are equicharacteristic or ramified,
the Positivity Conjecture implies (SP-1).  It is proved 
by applying Gabber's methods to the Positivity Conjecture.

\begin{tthm} \label{thm:KR}
Let $(R,\m)$ be a regular local ring with prime ideals $\p$ and 
$\q$ such that $\sqrt{\p+\q}=\m$ and $\dim(R/\p)+\dim(R/\q)=\dim(R)$.  If
$R$ is ramified or contains a field and $\chi(R/\p,R/\q)>0$, 
then $\p^{(m)}\cap\q\subseteq\m^{m+1}$ for all $m\geq 1$.
\end{tthm}

This theorem shows that (SP-1) holds in equal characteristic,
because the same is true of the Positivity Conjecture.
Because Kurano and Roberts expect the Positivity Conjecture to hold,
this theorem motivated them to formulate their conjecture, (SP-1).
It is worth noting that Kurano and Roberts~\cite[3.4]{kurano:pimsppi}
verify the conjecture directly for rings containing a field of 
characteristic 0, with no reference 
to positivity.  

The current author~\cite[2.3]{sather:isppi} has verified 
(SP-1) and (SP-2) when 
$R$ contains a field of arbitrary characteristic, also
with no reference to positivity.  This is accomplished 
by first proving that (ID-2) holds for rings containing a 
field, c.f., \cite[2.2]{sather:isppi}.  Then one verifies that (ID-2) implies
(SP-2), as follows.  If $R$ is a regular local ring, $\p$ a prime ideal 
of $R$ and $f$ a nonzero element of $\p$, then $e(R_{\p}/(f))=m$ if and 
only if $f\in\p^{(m)}\smallsetminus\p^{(m+1)}$.  If $\p$ and $\q$ are as 
in the statement of (SP-2)
and $f\in\p^{(m)}\cap\q^{(n)}\smallsetminus\m^{m+n}$, then we may apply 
(ID-2) to the hypersurface $R/(f)$ to arrive at a 
contradiction.  

The evidence in support of (ID-$i$) is the fact that it holds in 
equal characteristic.  It is worth noting that (ID-$i$) generalizes 
Serre's Intersection Theorem because, 
if $A$ is regular, then $e(A)=e(A_P)=e(A_Q)=1$.  One can take this 
as support for the conjecture as well.

\subsection{Relations between the conjectures} \label{subsec:imply}

In this subsection we discuss the implications in 
diagram (\ref{diagram}).  We have already seen how
(ID-2) implies (SP-2), and
a similar argument shows that (ID-1) implies (SP-1).
Clearly, (SP-2) implies (SP-1).  It is not immediately clear that (ID-2)
implies (ID-1), because $A/Q$ need not be analytically unramified in (ID-1).
However, (ID-2) implies (ID-1) for complete domains, and
the proof of our main theorem shows that we can reduce (ID-1) 
to the case of a complete domain.  It follows that (ID-2) implies (ID-1).

At this time, we do not know if the converses of the vertical
implications in diagram (\ref{diagram}) 
hold.  With Theorem~\ref{thm:KR} in mind, 
we remark that we do not know whether Serre's Positivity Conjecture implies
any of our conjectures, except (SP-1) in the ramified case.
Also, we do not know if any of our conjectures implies Positivity.

\subsection{Definitions and tools} \label{subsec:tools}

In this subsection, we give the
definition of the Hilbert-Samuel multiplicity and
catalog the properties we will use in the proof of our results. 

Let $(A,\n)$ be a local ring, $M$  a 
nonzero, finitely generated $A$-module of dimension $d$, 
and $I$ is an ideal of $A$ such that $M/IM$ 
has finite length.  The Hilbert function $P(n)=\len_A(M/I^{n+1}M)$ 
agrees with
a polynomial $Q(n)\in\mathbb{Q}[n]$
of degree $d$, for $n\gg 0$.  The Hilbert-Samuel multiplicity 
is the positive integer $e_A(I,M)$ such that
$Q(n)=\frac{1}{d!}e_A(I,M)n^d+(\text{lower degree terms})$.  When there is no 
ambiguity, we shall write $e(I,M)$.  If $I=\n$, we write 
$e_A(M)$ or $e(M)$.

The following is a list of properties satisfied by the Hilbert-Samuel 
multiplicity that are needed for the proof of our main theorem.
Let $A$ be a ring and $M$ a finite $A$-module.

\begin{property}  
\label{prop:1}
\textit{(Additivity Formula)}  Let $A$ be a local ring.  Then
\[ e_A(M)=\sum_{P}\len(M_P)e(A/P) \]
where the sum is taken over all primes $P$ of $A$ such that 
$\dim(A/P)=\dim(M)$.  
This sum is finite
because we need only take the sum over such primes which are also in the 
support of $M$.  (C.f., Bruns-Herzog~\cite[4.7.8]{bruns:cmr}.)
\end{property}

\begin{property}  
\label{prop:2}
Let $(R,\m)$ be a local subring of $A$ 
such that the extension $R\to A$ is module-finite.
Then $A$ is a semilocal ring such that $\dim(A)=\dim(R)$. 
Let $\{\n_1,\ldots,\n_n\}$ be the set of maximal ideals of $A$ 
such that $\n_i\cap R=\m$ and $\Ht(\n_i)=\dim(A)$.   Then
\[ e_{R}(\m,A)=
  \sum_i [A/\n_i:R/\m] e_{A_{\n_i}}\!(\m A_{\n_i},A_{\n_i}). 
    \]
(C.f., Nagata~\cite[(23.1)]{nagata:lr}.)
\end{property}

\begin{property}  
\label{prop:3}
Let $R$ be a subring of $A$ 
such that the extension $R\to A$ is module-finite, and 
let $\s$ be a prime ideal 
of $R$.  Then there are finitely many prime ideals of $A$ that
contract to $\s$ in $R$.
Let $\{ S_1,\ldots,S_j\}$ be the set of prime ideals of $A$ such that 
$S_i\cap R=\s$ and $\Ht(S_i)=\Ht(\s)$.  Then
\[ e_{R_{\s}}(\s R_{\s},A_{\s})=
  \sum_{i=1}^j [\kappa(S_i):\kappa(\s)] e_{A_{S_i}}\!(\s A_{S_i},A_{S_i})
    \]
where $\kappa(\s)$ and $\kappa(S_i)$ are 
the residue fields of $R_{\s}$ and $A_{S_i}$, respectively.
(Apply~\ref{prop:2} to the extension $R_{\s}\rightarrow A_{\s}$.)
\end{property}

\begin{property}  
\label{prop:4}
Let $A\rightarrow A'$ be a flat, local 
homomorphism of local rings $(A,\n)$ 
and $(A',\n')$ such that $\n A'=\n'$.  Then $e(A')=e(A)$.
(C.f., Herzog~\cite[2.3]{herzog:omlr}.)
\end{property}

\begin{property}  
\label{prop:5}
Let $A$ be a local ring and 
$I$ an ideal of $A$ such that $A/I$ has finite length.  If $M$ is an 
$A$-module of positive rank $r$, then 
$e_A(I,M)=e_A(I,A)\cdot r$.
(C.f., \cite[4.7.9]{bruns:cmr})
\end{property}

\begin{property}  
\label{prop:6}
If $A$ is local with infinite residue field, then 
there is an ideal $I$ generated by a system of parameters for $A$ 
such that $e(I,A)=e(A)$.
(C.f., Northcott-Rees~\cite[Theorem 1 of \S 6]{northcott:rilr}.)
\end{property}

\begin{property}  
\label{prop:7}
Assume that $(A,n)$ is local and equidimensional, and contains an
excellent local domain $(R,\m)$ such that the extension 
$R\rightarrow A$ is module-finite.  Let $P$ and $Q$ be prime ideals of $A$ 
such that $\sqrt{P+Q}=\n$,   
and let $\p=P\cap R$ and $\q=Q\cap R$.  If 
$e_R(\m,A)<e(A_P)+e(A_Q)$, then $\sqrt{\p+\q}=\m$.  
(C.f., \cite[2.1]{sather:isppi}.)
\end{property}

\subsection{Technical assumptions in the conjectures} \label{sec:ta}

  From the proof of our main theorem below, one can see why we make 
certain assumptions in (SP-$i$) and (ID-$i$).  Some of 
these assumptions are necessary;  we are not sure of the necessity of 
others.
In this subsection, we give some explanation and historical precedence 
for some of these assumptions.

Straightforward examples show that the following assumptions 
in (SP-$i$) are necessary:  (i) $R$ is regular, (ii) 
$\sqrt{\p+\q}=\m$, and (iii) $\dim(R/\p)+\dim(R/\q)=\dim(R)$.  
When we assume (i) and (ii), Serre's Intersection Theorem implies 
that $\dim(R/\p)+\dim(R/\q)\leq\dim(R)$.  That is, in (iii)
we are requiring the 
sum $\dim(R/\p)+\dim(R/\q)$ to achieve its maximal value.
Similar examples show that  we must 
assume that $\sqrt{P+Q}=\n$ in (ID-$i$).

In (ID-$i$), we assume that $A$ is quasi-unmixed.
There are 
straightforward examples of complete local rings that are not equidimensional 
for which the (ID-$i$) both fail, so the requirement 
that $A$ be equidimensional is certainly necessary.  
Since our arguments involve passing to the 
completion $A^*$, we must assume that $A^*$ is equidimensional.  
We do not know 
whether the full strength of the quasi-unmixedness assumption is needed.  
However, if one is inclined to assume that one's rings are excellent, then 
``quasi-unmixed'' is equivalent to ``equidimensional'', as we noted 
above, and is therefore necessary.

In (ID-$i$) we assume that $A/P$ 
is analytically unramified.  
The purpose of this assumption is to guarantee that 
the multiplicity $e(A_P)$ is well-behaved under passing to the 
completion $A^*$.  
To place this in an historical context, note that our proof is
very similar to Nagata's proof of the following result.

\begin{tthm} \label{thm:nagata} \rm (\cite[(40.1)]{nagata:lr})
\it Let $P$ be a prime ideal of a local ring $A$.  If $A/P$ is analytically 
unramified and if $\Ht(P)+\dim(A/P)=\dim(A)$, then $e(A_P)\leq e(A)$.
\end{tthm}

Regarding the analytically unramified assumption in this theorem,
Nagata (\cite[A2]{nagata:lr}) writes the following:
``It is not yet known to the writer's knowledge 
whether or not (40.1) is true without assuming that $P$ is analytically 
unramified.''  
If it is shown that this condition can be omitted from the 
statement of Theorem~\ref{thm:nagata}, then one will probably
be able to show that
the corresponding conditions 
should be omitted from the statements (ID-$i$).  

\section{Main Results} \label{sec:main}

The following theorem includes the Main Theorem announced in 
Section~\ref{sec:intro}.

\begin{thm} \label{thm:reduce}
With the notation of Conjectures~\ref{conj:sp} and~\ref{conj:id}, 
fix $i=1$ or $2$.  Then 
the following conditions are equivalent.
\begin{enumerate}
\item \emph{(SP-$i$)} holds for all regular local rings. \label{item:1}
\item \emph{(SP-$i$)} holds for all complete, 
unramified regular local rings of 
mixed characteristic.  \label{item:2}
\item For every prime number $p$, every complete $p$-ring 
$(V, pV)$ and every $d\geq 1$, \emph{(SP-$i$)} holds for 
$V[X_1,\ldots,X_d]_{(p,X_1,\ldots,X_d)}$. \label{item:3}
\item \emph{(ID-$i$)} holds for all quasi-unmixed local rings. \label{item:4}
\item \emph{(ID-$i$)} holds for all hypersurfaces 
over an arbitrary complete, unramified 
regular local ring of mixed characteristic. \label{item:5}
\item For every prime number $p$, every complete $p$-ring 
$(V, pV)$ and every $d\geq 1$, \emph{(ID-$i$)} holds 
for all hypersurfaces over 
$V[X_1,\ldots,X_d]_{(p,X_1,\ldots,X_d)}$.    \label{item:6}
\end{enumerate}
\end{thm}

In the statement, a \textit{$p$-ring} is a discrete valuation ring $V$ 
whose maximal ideal is generated by the prime number $p$.


\begin{proof}
We summarize the implications we shall prove in the following diagram.
\[
\xymatrix{
(\ref{item:1}) \ar@{=>}[r] & (\ref{item:2}) \ar@{<=>}[r] \ar@{<=>}[d] 
  & (\ref{item:3}) \ar@{<=>}[d] \\
(\ref{item:4}) \ar@{=>}[u] \ar@{<=>}[r] & (\ref{item:5}) \ar@{=>}[r] 
  & (\ref{item:6}) } 
\]


The implications ``(\ref{item:1}) $\implies$ (\ref{item:2})'' and 
``(\ref{item:4})$\implies$(\ref{item:5})'' are obvious.  The implication
``(\ref{item:2})$\implies$(\ref{item:3})''  follows by passing to the 
completion $R^*$ of $R$;  use the faithful flatness of $R\rightarrow R^*$ 
as in the first step of ``(\ref{item:3})$\implies$(\ref{item:2})''.  
The implication
``(\ref{item:5})$\implies$(\ref{item:6})'' also follows by passing to the 
completion;  argue as in the proof of~\cite[2.3]{sather:isppi}.
As noted in the introduction, if $R$ is a regular local ring with prime 
ideal $\p$ and 
$f$ is a nonzero element of $\p$, then 
$f\in\p^{(m)}\smallsetminus\p^{(m+1)}$ if and only if $e(R_{\p}/(f))=m$.  
  From this the implications ``\ref{item:4} $\implies$ \ref{item:1}'',
``(\ref{item:2})$\iff$(\ref{item:5})'' and 
``(\ref{item:3})$\iff$(\ref{item:6})'' follow easily.

``(\ref{item:5})$\implies$(\ref{item:4})''.  
We verify this implication in 
the case $i=2$;  the case $i=1$ is similar.
Let $(A,\n)$ be a quasi-unmixed, local ring with prime ideals $P$ 
and $Q$ such 
that $A/P$ and $A/Q$ are analytically unramified, $\sqrt{P+Q}=\n$, and 
$e(A)<e(A_P)+e(A_Q)$.  


We now show that we may assume without loss of generality that $A$ is 
a complete local domain of mixed characteristic with infinite residue field.


Step 1.  Pass to the ring $A(X)=A[X]_{\n A[X]}$ to assume that the residue 
field of $A$ is infinite.  Let $\n(X)=\n A(X)$, $P(X)=PA(X)$ and $Q(X)=QA(X)$.  
Then $A(X)$ is a local ring with maximal ideal $\n(X)$.  
Multiplicities are preserved by~\ref{prop:4}.
That is, 
$e(A(X))=e(A)$, $e(A(X)_{P(X)})=e(A_P)$ and $e(A(X)_{Q(X)})=e(A_Q)$.
Both rings $A(X)/P(X)$ 
and $A(X)/Q(X)$ are analytically unramified by~\cite[(36.8)]{nagata:lr}. 
Also, $A(X)$ is quasi-unmixed:  this is equivalent to it being 
universally catenary and equidimensional by~\cite{ratliff:qldafccpiII},
and since $A$ satisfies these properties, the same is true of
$A(X)$.  If the result holds in $A(X)$, then
\begin{align*}
\dim(A/P)+\dim(A/Q)&=\dim(A(X)/P(X))+\dim(A(X)/Q(X))\\
&\leq\dim(A(X))=\dim(A) 
\end{align*}
as desired.  

Step 2.  Pass to the completion $(A^*,\n^*)$ to assume that $A$ is complete 
and equidimensional with infinite residue field.  
Let $P^*$ be a prime ideal of $A^*$ that is minimal over $PA^*$ such 
that $\Ht(P^*)=\Ht(P)$, and similarly for $Q^*$.  Since $A/P$ is 
analytically unramified, $PA^*_{P^*}=P^*A^*_{P^*}$.  The fact that the 
extension $A_P\rightarrow A^*_{P^*}$ is flat therefore implies that 
$e(A^*_{P^*})=e(A_P)$ by~\ref{prop:4}.  
Similarly, $e(A^*_{Q^*})=e(A_Q)$.  If the result 
holds for $A^*$, then it holds for $A$, as in the previous step.

Step 3.  Pass to 
the quotient $A/I$ for a suitably chosen minimal prime $I$, 
to assume that $A$ is a complete domain with 
infinite residue field.  To make this 
reduction, it suffices to verify that 
there is a minimal prime $I$ of $A$ contained in 
$P\cap Q$ such that $e(A/I)<e(A_P/IA_P)+e(A_Q/IA_Q)$.  

First, we show how this gives the desired reduction.  Let $A'=A/I$ with 
prime ideals $\n'=\n A'$, $P'=PA'$ and $Q'=\q A'$.  Since $A$ is 
equidimensional, $\dim(A')=\dim(A)$, and since $I\subseteq P\cap Q$, it 
follows that $\dim(A'/P')=\dim(A/P)$ and $\dim(A'/Q')=\dim(A/Q)$.  
Therefore, as in Step 1, we may pass to $A'$.

Now we prove that such a minimal prime $I$ exists.  
Let $\{I_1,\ldots,I_g\}=\min(A)$.  
Suppose that $e(A/I_j)\geq e(A_P/I_jA_P)+e(A_Q/I_jA_Q)$ for every $j$ such 
that $I_j\subseteq P\cap Q$.  (This supposition includes 
the hypothetical possibility that no $I_j$ is contained in $P\cap Q$.)
By~\cite[(40.1)]{nagata:lr}, $e(A_P/I_jA_P)\leq e(A/I_j)$ for every $I_j$ 
contained in $P$, and similarly for $I_j$ contained in $Q$.
The Additivity Formula then implies that
\begin{align*}
e(A_P)+e(A_Q)
  &=\sum_{I_j\subseteq P}e(A_P/I_jA_P)\len(A_{I_j})
    +\sum_{I_j\subseteq Q}e(A_Q/I_jA_Q)\len(A_{I_j}) \\
  &=\sum_{\substack{I_j\subseteq P \\ I_j\not\subseteq Q}}
          e(A_P/I_jA_P)\len(A_{I_j})  
   +\sum_{\substack{I_j\subseteq Q \\ I_j\not\subseteq P}}
          e(A_Q/I_jA_Q)\len(A_{I_j})  \\
  &\quad +\sum_{I_j\subseteq P\cap Q} 
         [e(A_P/I_jA_P)+e(A_Q/I_jA_Q)]\len(A_{I_j}) \\
  &\leq\sum_{\substack{I_j\subseteq P \\ I_j\not\subseteq Q}}
          e(A/I_j)\len(A_{I_j})  
   +\sum_{\substack{I_j\subseteq Q \\ I_j\not\subseteq P}}
          e(A/I_j)\len(A_{I_j})  \\
  &\quad +\sum_{I_j\subseteq P\cap Q} 
         e(A/I_j)\len(A_{I_j}) \\
  &\leq\sum_j e(A/I_j)\len(A_{I_j}) \\
  &=e(A).
\end{align*}
This clearly contradicts the 
assumption that $e(A)<e(A_P)+e(A_Q)$.  Thus, 
there is a minimal prime $I$ of $A$ such that 
$I\subseteq P\cap Q$ and $e(A/I)<e(A_P/IA_P)+e(A_Q/IA_Q)$, as claimed.

Since (ID-$i$) has been verified for rings containing a field, we 
assume without loss of generality that $A$ has characteristic 0 and that 
$k=A/\n$ has characteristic $p>0$.  

%
%

Let $V$ be a coefficient ring 
for $A$, which is a complete $p$-ring.  Since the residue field of $A$ is 
infinite, \ref{prop:6} implies that there exists
system of parameters $\x=x_1,\ldots,x_d$ of 
$A$ such that the ideal $J=\x A$ satisfies
$e(J, A)=e(A)$.  Let $R=V[\![x_1,\ldots,x_d]\!]$ which is a complete domain contained 
in $A$ such that the extension $R\rightarrow A$ is module finite, and 
the induced map on residue fields is an isomorphism.
Let $\m=(p,x_1,\ldots,x_d)R=\n\cap R$, $\p=P\cap R$ and $\q=Q\cap R$.
Since $J=\x A\subseteq\m A\subseteq\n$, we have 
$e(A)=e(J,A)\geq e(\m A,A)\geq e(\n,A)=e(A)$.  
Therefore,
\[ e(A)=e_A(\m A,A)[A/\n,R/\m]=e_R(\m,A)=\rank_R(A) e(R) \]
by~\ref{prop:2}.
Since $e(A)<e(A_P)+e(A_Q)$, 
\ref{prop:7} implies that $\sqrt{\p+\q}=\m$.

If $R$ is a regular local ring, then Serre's Intersection Theorem implies 
that
\[ \dim(A/P)+\dim(A/Q)=\dim(R/\p)+\dim(R/\q)\leq\dim(R)=\dim(A) \]
as desired.  Therefore, we assume 
that $R$ is not regular.  Let $X_1,\ldots,X_d$ be indeterminates and let 
$B=V[\![X_1,\ldots,X_d]\!]$.  Then $B$ surjects onto $R$ via the natural 
homomorphism $\pi$ which sends $X_i$ to $x_i$.  
Since $\dim(R)=d=\dim(B)-1$ and $R$ is a domain, the kernel of $\pi$ is a 
height 1 prime of $B$ and is therefore principal.  Thus, $R$ is a 
hypersurface over $B$.  If we can show that $e(R)<e(R_{\p})+e(R_{\q})$, 
then our assumptions imply that $\dim(R/\p)+\dim(R/\q)\leq \dim(R)$, 
completing the proof.

Let $r=\rank_R(A)=\rank_{R_{\p}}(A_{\p})>0$ 
so that $e(A)=re(R)$.  We claim that 
$re(R_{\p})\geq e(A_P)$.  As before,
\begin{align*}
re(R_{\p}) & =\rank_{R_{\p}}(A_{\p}) e(\p R_{\p},R_{\p}) \\
            & =e(\p R_{\p},A_{\p}) \\
  & =\sum_{P_i\cap R=\p} 
         e_{A_{P_i}}\!(\p A_{P_i},A_{P_i})[\kappa(P_i):\kappa(\p)] \\
\intertext{(where the sum is taken over all primes $P_i$ of $A$ contracting 
to $\p$ in $R$)}
  & \geq e_{A_P}(\p A_P,A_P) \\
  & \geq e_{A_P}(PA_P,A_P)\\
  & =e(A_P).
\end{align*}
Similarly, $re(R_{\q})\geq e(A_Q)$.  
Thus, 
\[ re(R)=e(A)<e(A_P)+e(A_Q)\leq r(e(R_{\p})+e(R_{\q})) \]
and since $r>0$, we have $e(R)<e(R_{\p})+e(R_{\q})$, as desired.

``(\ref{item:3})$\implies$(\ref{item:2})''.  The main ideas for this proof 
come from Hochster~\cite{hochster:tithtomocr} and 
Dutta~\cite{dutta:tsbqcgims}, each of which relies heavily on Artin 
approximation~\cite{artin:aasclr} and on the work of 
Peskine-Szpiro~\cite{peskine:dpfcl}. 

Let $V$ be a complete $p$-ring and $d\geq 1$, and suppose that there were a 
counterexample to (SP-2) in $R=V[\![X_1,\ldots,X_d]\!]$.  
That is, there exist
prime ideals $\p,\q$ in $R$ and an element $f\in R$ such that 
$\sqrt{\p+\q}=\m$, 
$\dim(R/\p)+\dim(R/q)=d+1$ and 
$f\in\p^{(m)}\cap\q^{(n)}\smallsetminus\m^{m+n}$.  
Since our 
conjectures are easily verified for rings of dimension less than 2, 
we have $\dim(R)\geq 2$.

We may assume without loss of generality that $k=V/pV$ is 
algebraically closed, as follows.  Let $k'$ denote the algebraic closure 
of $k$.  By Matsumura~\cite[29.1]{matsumura:crt}
there exists a $p$-ring $V'$ such that $V'/pV'\cong k'$ and a flat ring 
homomorphism $V\to V'$.  Replace $V'$ by its completion to assume 
that $V'$ is a complete $p$-ring.  Let $R'=V'[\![X_1,\ldots,X_d]\!]$, which 
is a regular local ring of dimension $d+1=\dim(R)$
with 
maximal ideal $\m'=(p,X_1,\ldots,X_d)R'$.
The induced ring homomorphism $R\to R'$ is flat and local, so we may
fix prime ideals $\p',\q'\subset R'$ such that 
$\p'\cap R=\p$ and $\q'\cap R=\p$.
It follows that $\dim(R'/\p')=\dim(R/\p)$ and $\dim(R'/\q')=\dim(R/\q)$,
since going-down holds between $R$ and $R'$, and each ring is a 
catenary domain;  apply~\cite[15.1]{matsumura:crt}.
For $i=1,\ldots,d$ and $q\gg 0$ we have
$p^q,X_i^q\in\p+\q\subseteq\p'+\q'$, and it follows that 
$\sqrt{\p'+\q'}=\m'$.
Since $f\in\p^{(m)}$, 
there exists 
$s\in R\smallsetminus\p$ such that $sf\in\p^m$.  Since $\p'\cap R=\p$,
it follows that
$s\not\in\p'$, and the fact that $sf\in\p^m\subseteq(\p')^m$ implies that 
$f\in(\p')^{(m)}$.  Similarly, $f\in(\q')^{(n)}$.  Finally, if 
$f\in(\m')^{m+n}$, then $f\in(\m')^{m+n}\cap R=\m^{m+n}$ by faithful 
flatness, a contradiction.  Thus, the counterexample in $R$ ascends to a 
counterexample in $R'$, and we may assume that $V/pV$ is algebraically 
closed.

Let 
$R_0=V[X_1,\ldots,X_d]_{(p,X_1,\ldots,X_d)}$ with maximal ideal $\m_0$.  
We shall show that the counterexample in $R$ gives rise to a 
counterexample in $R_0$.  
More specifically, we shall show that there 
exist prime ideals $\p_0,\q_0$ in $R_0$ and $g\in R_0$ such that 
$\sqrt{\p_0+\q_0}=\m_0$, $\dim(R_0/\p_0)=\dim(R/\p)$, 
$\dim(R_0/\q_0)=\dim(R/\q)$, and 
$g\in\p_0^{(m)}\cap\q_0^{(n)}\smallsetminus\m_0^{m+n}$.  

First, we show how the given counterexample gives rise to 
a counterexample in the Henselization $R_0^h$ of $R_0$, and therefore in a 
pointed \'{e}tale neighborhood $R_1$ of $R_0$.  
This is accomplished by 
showing that the essential data describing the counterexample is 
determined by a finite number of polynomial equations (necessarily in a
finite number of variables) with coefficients in $R_0$.  Artin's 
approximation theorem will then yield a solution to these equations in 
$R_0^h$.  Furthermore, given an integer $e\geq 1$ 
we can guarantee that 
the solutions in $R_0^h$ will be congruent to the original solutions 
modulo $\m^e$.
Since there are finitely many variables, the solution will lie in 
a pointed \'{e}tale neighborhood $R_1$ of $R_0$.  
Note that $R_0^h$ and 
each pointed 
\'{e}tale neighborhood $R_1$ of $R_0$ is an unramified regular local ring of 
dimension $d+1$ with regular system of parameters $p,X_1,\ldots,X_d$.
Let 
$y_1,\ldots,y_a\in\p$ be a minimal set of generators of $\p$ and 
$z_1,\ldots,z_b\in\q$ a minimal set of generators of $\q$.  

1.  A finite number of equations determine the 
minimal free resolutions of $R/\p$ and $R/\q$;  
c.f., \cite[6.2]{peskine:dpfcl}.  Let 
$Y_1,\ldots,Y_a,Z_1,\ldots,Z_b$ be the variables that ``keep track'' of the 
elements $y_1,\ldots,y_a,z_1,\ldots,z_b$.  Applying Artin approximation 
will yield ideals $I, J$ in a pointed \'{e}tale neighborhood $(R_1,\m_1)$ 
of $R_0$ that are generated by the solutions in $R_0$ that are substituted 
for the variables $Y_1,\ldots,Y_a,Z_1,\ldots,Z_b$.

2.  A finite number of equations determine the fact 
that $\sqrt{\p+\q}=\m$.  There exists $q\geq 1$ such that 
$p^q,X_j^q\in\p+\q$ for $j=1,\ldots,d$, so we may write 
\[ p^q=\sum_{i=1}^a c_iy_i+\sum_{i=1}^bc_i'z_i
\qquad\text{and}\qquad 
X_j^q=\sum_{i=1}^a c_{j,i}y_i+\sum_{i=1}^bc_{j,i}'z_i \]
for $j=1,\ldots,d$ and 
$c_i, c_{j,i},c'_i, c'_{j,i}\in R$.  Using variables $C_i, C_{j,i},
C'_i, C'_{j,i}$ for the 
coefficients, and the variables $Y_i, Z_i$ already designated, we can 
express this information in $d+1$ equations.

3.  A finite number of equations determine the dimension 
of $R/\p$ and $R/\q$.  This is a result of 
Hochster~\cite{hochster:tithtomocr}.  A proof can be found 
in~\cite[8.4.4]{bruns:cmr}.  See~\cite[3.5]{dutta:tsbqcgims} for a 
different proof.

4.  A finite number of equations determine the fact that 
a given sequence $\alpha_1,\ldots,\alpha_r$ is a system of parameters for 
$R/\p$.  We can keep track of $\dim(R/\p)$.  We can keep track of the 
fact that the $\alpha_i$ are elements of $\m$ and that the
sequence has $\dim(R/\p)$ elements.  Also we can keep track of 
the fact that $(R/\p)/(\alpha_1,\ldots,\alpha_r)$ has dimension 0.  

It 
follows that we can keep track of whether a given element $\alpha_1\in\m$ is 
part of a system of parameters for $R/\p$ by extending it to a full system 
of parameters.

5.  A finite number of equations determine the fact that 
$f\in\p^{(m)}$.  There are two cases.  

Case 1.  If $f\in\p^{(m)}\smallsetminus\p^m$, then there exists 
$s\in\m\smallsetminus\p$ such that $sf\in\p^m$.  That is, we can write
\[ sf=\!\!\!\!\sum_{i_1,\ldots,i_a\geq 0}\!\!\!\!
c''_{i_1,\ldots,i_a}y_1^{i_1}\cdots y_a^{i_a} \]
where the sum is finite and the coefficients are in $R$.  Since 
$s\in\m\smallsetminus\p$, $s$ is part of a system 
of parameters for $R/\p$;  we have already noted that we can keep track of 
this fact with a finite number of equations.  By using variables 
$C''_{i_1,\ldots,i_a}$ we can also keep track of the fact that 
$sf\in\p^{m}$.  Let $S,F$ be the variables we use to keep track of the 
elements $s,f$.  

Case 2.
If $f\in\p^m$, then we can keep track of the fact 
that $f\in\p^m$ with a single equation, 
as in case 1.  Using $s=1$ and the equation $S=1$ 
keeps this case consistent with case 1.

In our application of Artin approximation, we may then conclude that there 
are elements $s_1,f_1\in R_1$ such that $s_1f_1\in I^m$ and either $s_1=1$ or
$s_1$ is part of a system of parameters for $R_1/I$.  
Similarly, there 
is an element $t_1\in R_1$ such that $t_1f_1\in J^n$ and either $t_1=1$ or
$t_1$ is part of a system of parameters for $R_1/J$.  

6.  Choosing $e\geq m+n$, we require that the solutions in $R_1$ are 
congruent to the original solutions in $R$ modulo $\m^{m+n}$.  In 
particular, since $f\not\in\m^{m+n}$, it follows that $f_1\not\in\m^{m+n}$ 
and, therefore, that $f_1\not\in\m_1^{m+n}$.

To summarize, we have shown the existence of the following:
\begin{itemize}
\item[(a)]  a pointed \'{e}tale neighborhood $(R_1,\m_1)$ 
of $(R_0,\m_0)$ and an element $f_1\in 
R_1\smallsetminus\m_1^{m+n}$;
\item[(b)]  ideals $I,J\subset R_1$ such that  
$\sqrt{I+J}=\m_1$, $\dim(R_1/I)=\dim(R/\p)$ and $\dim(R_1/J)=\dim(R/\q)$;
\item[(c)]  an element $s_1\in R_1$ such that $s_1f_1\in I^m$ and either $s_1=1$ or
$s_1$ is part of a system of parameters for $R_1/I$;  
\item[(d)]  an element $t_1\in R_1$ such that $t_1f_1\in J^n$ and either $t_1=1$ or
$t_1$ is part of a system of parameters for $R_1/J$.  
\end{itemize}
In particular, $R_1$ is a regular local 
ring of dimension $d+1$, essentially of finite type and smooth over $V$,
$\m_1=\m_0 R_1$, and $R_1/\m_1=V/pV$.  Let $\p_1$ be a minimal 
prime of $R_1/I$ such that $\dim(R_1/\p_1)=\dim(R_1/I)$, and let $\q_1$ 
be a minimal prime of $R_1/J$ such that $\dim(R_1/\q_1)=\dim(R_1/J)$.  
Then $\sqrt{\p_1+\q_1}=\m_1$ and $\dim(R_1/\p_1)+\dim(R_1/\q_1)=\dim(R_1)$.
Since $s_1$ is either part of a system of parameters on $R_1/I$ or $s=1$, 
it follows that $s_1\not\in\p_1$.  Since $s_1f_1\in I^m\subseteq\p^m$, it 
follows that $f_1\in\p^{(m)}$.  Similarly, $f_1\in\q_1^{(n)}$.  Thus, 
we have a counterexample in $R_1$, as desired.  

Next, we show that
we may assume without loss of 
generality that $f_1$ is irreducible in our conjecture.  Since 
$R$ is a unique factorization domain, write $f_1=F_1\cdots F_v$ with each
$F_i\in\m_1$ irreducible.  
For $i=1,\ldots,v$ let $m_i,n_i\geq 0$ be the unique integers such that 
$F_i\in\p_1^{(m_i)}\smallsetminus\p_1^{(m_i+1)}$ and 
$F_i\in\q_1^{(n_i)}\smallsetminus\q_1^{(n_i+1)}$.  In other words, let 
$m_i=e((R_1)_{\p_1}/(F_i))$ and $n_i=e((R_1)_{\q_1}/(F_i))$.  
Assuming that 
$f_1\in\p_1^{(m)}\cap\q_1^{(n)}\smallsetminus\m^{m+n}$, we claim that there 
exists an integer $i$ such that $1\leq i\leq v$ and 
$F_i\not\in\m_1^{m_i+n_i}$.  This will show that our counterexample 
gives rise to an irreducible counterexample.  Suppose that
each $F_i\in\m_1^{m_i+n_i}$.  
(If $m_i=0$ or $n_i=0$, this is 
automatic by~\cite[(38.3)]{nagata:lr}.)  From this it follows that
\[ f_1=F_1\cdots F_v\in\m_1^{\sum_i (m_i+n_i)}.\]
The Additivity Formula implies 
that
\[ e((R_1)_{\p_1}/(f_1))=\sum_i e((R_1)_{\p_1}/(F_i))=\sum_i m_i \]
or, in other words,
\[ f_1\in\p_1^{(\sum_i m_i)}\smallsetminus\p_1^{(1+\sum_i m_i)}. \]
Since $f_1\in\p_1^{(m)}$, it follows that $m\leq\sum_i m_i$.  Similarly, 
$n\leq\sum_i n_i$.  Therefore, $m+n\leq \sum_i(m_i+n_i)$, and it follows 
that 
\[ f_1\in\m_1^{\sum_i (m_i+n_i)}\subseteq\m_1^{m+n} \]
which contradicts our hypothesis.  

We also note that
$f_1\in\m_1^2$, as follows.  Suppose not.  Since 
$f_1\in\p_1\cap\q_1\subseteq\m_1$, the fact that $f_1\not\in\m_1^2$ implies 
that $R'=R_1/(f_1)R_1$ is a regular local ring of dimension $d$ with prime 
ideals $\p'=\p_1/(f_1)R_1$ and $\q'=\q/(f_1)R_1$ such that 
$\sqrt{\p'+\q'}=\m'$, and 
\[ \dim(R'/\p')+\dim(R'/\q')=\dim(R_1/\p_1)+\dim(R_1/\q_1)=d+1>\dim(R'). \]
This contradicts Serre's Intersection Theorem.  
Observe that $f\not\in pR_1$ because
$f_1$ is irreducible and $f\in\m^2$.

Finally, we show how 
the following theorem of Dutta~\cite[1.3]{dutta:tsbqcgims} gives rise 
to a 
counterexample in $R_0$, thus completing our proof.  

\begin{thm} 
Let $(A,M,K)$ be a regular local ring of dimension $d+1$, essentially of 
finite type and smooth over an excellent discrete valuation ring
$(U,\pi)$ such that 
$K$ is separably generated over $U/\pi U$.  Let $a (\neq 0)\in M^2$ be 
such that $a\not\in\pi A$ (i.e., $\{\pi,a\}$ form an $A$-sequence).  Then 
there exists a regular local ring $(B,N,K)\subset (A,M,K)$ such that

1.  The ring
$B$ is a localization of a polynomial ring $W[Y_1,\ldots,Y_d]$ at a 
maximal ideal of the type $(\pi,h(Y_1),Y_2,\ldots,Y_d)$ where $h$ is a 
monic irreducible polynomial in $W[Y_1]$ and $(W,\pi)$ is an excellent 
discrete valuation ring contained in
$A$;  moreover, $A$ is an \'{e}tale neighborhood of $B$.

2.  There exists an element $g$ in $B\cap aA$ such that $B/gB\to A /aA$ 
is an isomorphism.  Furthermore $gA=aA$.
\end{thm}

We apply this
theorem to the ring $A=R_1$ and the element $a=f_1$.  
The discrete valuation ring $W$ whose existence is guaranteed by the 
theorem is constructed by choosing elements $y_1,\ldots,y_t\in A$ 
such that their residues modulo $M$ form 
a transcendence basis of $K$ over $U/\pi U$ and then setting
$W=V[y_1,\ldots,y_t]_{\pi V[y_1,\ldots,y_t]}$.  
The fact that 
$V/pV\cong R_1/\m_1$ implies that we may take 
$W=V$ in the conclusion.  Since $V/pV$ is algebraically closed, 
the maximal ideal 
$(\pi,h(Y_1),Y_2,\ldots,Y_d)$ will be of the form 
$(\pi,Y_1-z_1,Y_2,\ldots,Y_d)$ for some element $z_1\in V$;  therefore, $B\cong 
R_0$.  Without loss of generality, we may assume that $B=R_0$.
Let $\p_0=R_0\cap \p_1$ and $\q_0=R_0\cap \q_1$.  Then 
$g\in\p_0\cap\q_0$.  Since the induced map $R_0/gR_0\to R_1/f_1R_1$ is an 
isomorphism, we see that $R_0/\p_0\cong R_1/\p_1$ and $R_0/\q_0\cong 
R_1/\q_1$.  Furthermore, under this isomorphism, 
$(\p_0/gR_0)+(\q_0/gR_0)\cong (\p_1/f_1R_1)+(\q_1/f_1R_1)$, which implies that
$\sqrt{\p_0+\q_0}=\sqrt{\p_0+\q_0+gR_0}=\m_0$.  Also, $(R_0/gR_0)_{\p_0}\cong 
(R_1/f_1R_1)_{\p_1}$, so that $e((R_0/gR_0)_{\p_0})=e((R_1/f_1R_1)_{\p_1})$;  
in 
particular, $g\in\p_0^{(m)}$.  Similarly, $g\in\q_0^{(n)}$ and 
$g\not\in\m_0^{m+n}$.  Thus, our counterexample in $R_1$ gives rise to a 
counterexample in $R_0$, completing the proof.
\end{proof}

Using methods similar to those employed in the proof of 
Theorem~\ref{thm:reduce}, we verify the following weaker version of 
(ID-2).

\begin{thm} \label{thm:weak(ID-2)}
Let $(A,\n)$ be a quasi-unmixed local ring with 
prime ideals $P$ and $Q$ such that $\sqrt{P+Q}=\n$.  If $e(A)<e(A_P)+e(A_Q)$,
e.g., if $e(A)=e(A_P)$, then $\dim(A/P)+\dim(A/Q)\leq\dim(A)+1$.
\end{thm}


\begin{proof}
Let $A, \n, P, Q$ be as in the statement of the theorem.  As in the proof 
of the implication ``(\ref{item:5})$\implies$(\ref{item:4})'' in 
Theorem~\ref{thm:reduce} we may assume that $A$ is a hypersurface 
over the ring $B=V[\![X_1,\ldots,X_d]\!]$, with surjection $\pi:B\to A$.  
Let $\mm=\pi^{-1}(\n)$ which is the maximal ideal of $B$, 
$\pp=\pi^{-1}(P)$ and $\qqq=\pi^{-1}(P)$.  Then  $\pp+\qqq$ is 
$\mm$-primary, so by Serre's Intersection 
Theorem, 
\[ \dim(A/P)+\dim(A/Q)=\dim(B/\pp)+\dim(B/\qqq)\leq\dim(B)=\dim(A)+1 \]
as desired.
\end{proof}





\section{An example} \label{sec:ex}

With Conjecture~\ref{conj:sp} in mind, it is natural to ask, ``What is the 
correct conjecture to make for symbolic powers in nonregular rings?''  This 
is a tricky business, as symbolic powers can behave badly 
in general.  For example, in the ring $k[\![X,Y]\!]/(XY)=k[\![x,y]\!]$,
the prime ideal $\p=(x)$ satisfies $\p^{(m)}=\p$ for all $m\geq 1$.  
However, a few interesting suggestions have been 
offered in this direction.  For instance, we have Question~\ref{q} 
below, 
which is motivated by the work of Hochster-Huneke~\cite{hochster:csopi}.
There, the authors use tight closure methods to investigate symbolic
powers.

For an ideal $I$ in a ring $R$ of prime characteristic $p>0$, let 
$I^*$ denote the tight closure of $I$.  A thorough introduction to 
the theory of tight closure can be found in~\cite[Chapter 10]{bruns:cmr}
the monograph of
Huneke~\cite{huneke:tcaia}.

\begin{question} \label{q}
Let $k$ be an algebraically closed field of prime characteristic $p>0$ and 
\[ R=k[X_1,\ldots,X_d]/(f)=k[x_1,\ldots,x_d] \]
for some                  
prime element $f\neq 0$ and $\m=(x_1,\ldots,x_d)R$.
Let $\p,\q$ prime ideals of $R$ such                         
that $\p+\q=\m$ and              
$\dim(R/\p)+\dim(R/\q)=\dim(R)+1$.  
Does the containment
\[ \p^{(m)}\cap\q\subseteq(\m^{m+1})^* \]
hold for all $n\geq 1$?
\end{question}

For $n\leq 3$, the answer to this question is ``yes''.
However, for $n=4$ the answer is ``no'', as the following 
example shows.

\begin{ex}
Fix an integer $s\geq 3$, and let 
\[ R=k[X,Y,Z,U]/(XY(Z+U)-U^sZ)=k[x,y,z,u]. \] 
Applying Eisenstein's criterion to $f=XY(Z+U)-U^sZ\in k[Y,Z,U][X]$ 
shows that $f$ is prime.
Let $\m=(x,y,z,u)R$, $\p=(x,u)R$ 
and $\q=(y,z)R$.  Then $\p+\q=\m$ and $\dim(R/\p)+\dim(R/\q)=4=\dim(R)+1$.  
By Hochster-Huneke~\cite[(6.2)]{hochster:fteasbc}, every element in the
Jacobian ideal has a power that is a test element.
In this example, the Jacobian ideal is 
\[ J=(y(u+z),x(u+z),xy-su^{s-1}z,xy-u^s)R.\] 
Fix an integer $q\geq 1$ such that $(xy-u^s)^q$ is a test element and 
$2q$ is divisible by $s-1$.  Set 
\[ m=\frac{2q}{s-1}+1>\frac{2q+1}{s-1}>0.\]
We claim that 
$x^my\in (\p^{(ms)}\cap \q)\smallsetminus (\m^{ms+1})^*$.

First, we show that $x^my\in \p^{(ms)}\cap \q$.  It 
suffices to show that $x^m\in\p^{(ms)}$.  This follows from the fact 
that $y^m(z+u)^m\not\in\p$ and $x^my^m(z+u)^m=u^{ms}z^m\in\p^{ms}$.

Next, we show that $x^my\not\in(\m^{ms+1})^*$.  As $(xy-u^s)^q$ is a 
test element, it suffices 
to show that $(xy-u^s)^qx^my\not\in\m^{ms+1}$.  
Expanding $(xy-u^s)^qx^my$ yields
\[ (xy-u^s)^qx^my=
     x^{q+m}y^{q+1}+[-qx^{q+m-1}y^{q}u^s+\cdots+(-1)^qx^myu^{sq}]. \]
The quantity in brackets is in $\m^{ms+1}$ by our choice of $m$.  It 
remains to show that $x^{q+m}y^{q+1}\not\in\m^{ms+1}$.  This is 
straightforward and can be verified by the interested reader.
\end{ex}


\noindent Acknowledgements.  
I am grateful to L.~Avramov, S.~Dutta, P.~Griffith, C.~Huneke, 
D.~Katz, D.~White and the referee for their 
helpful comments and suggestions.


\begin{thebibliography}{10}

\bibitem{artin:aasclr}
M.~Artin, \emph{Algebraic approximation of structures over complete local
  rings}, Inst. Hautes \'Etudes Sci. Publ. Math. (1969), no.~36, 23--58.

\bibitem{berthelot:ava}
P.~Berthelot, \emph{Alt\'erations de vari\'et\'es alg\'ebriques (d'apr\`es {A}.
  {J}. de {J}ong)}, Ast\'erisque (1997), no.~241, 273--311.

\bibitem{bruns:cmr}
W.~Bruns and J.~Herzog, \emph{Cohen-{M}acaulay rings}, Cambridge Studies in
  Advanced Mathematics, vol.~39, Cambridge University Press, Cambridge, 1993.

\bibitem{dutta:tsbqcgims}
S.~P. Dutta, \emph{A theorem on smoothness---{B}ass-{Q}uillen, {C}how groups
  and intersection multiplicity of {S}erre}, Trans. Amer. Math. Soc.
  \textbf{352} (2000), no.~4, 1635--1645.

\bibitem{eisenbud:essfi}
D.~Eisenbud and B.~Mazur, \emph{Evolutions, symbolic squares, and {F}itting
  ideals}, J. Reine Angew. Math. \textbf{488} (1997), 189--201.

\bibitem{gillet:knmi}
H.~Gillet and C.~Soul{\'e}, \emph{{$K$}-th\'eorie et nullit\'e des
  multiplicit\'es d'intersection}, C. R. Acad. Sci. Paris S\'er. I Math.
  \textbf{300} (1985), no.~3, 71--74.

\bibitem{herzog:omlr}
B.~Herzog, \emph{On the {M}acaulayfication of local rings}, J. Algebra
  \textbf{67} (1980), 305--317.

\bibitem{hochster:nimrrlr}
M.~Hochster, \emph{Nonnegativity of intersection multiplicities in ramified
  regular local rings following {G}abber/{D}e {J}ong/{B}erthelot}, unpublished
  notes.

\bibitem{hochster:tithtomocr}
M.~Hochster, \emph{Topics in the homological theory of modules over commutative
  rings},  CBMS Conference Ser.~in Math., No.~24,
  Amer.~Math.~Soc., Providence, R.I., 1975.

\bibitem{hochster:fteasbc}
M.~Hochster and C.~Huneke, \emph{{$F$}-regularity, test elements, and smooth
  base change}, Trans. Amer. Math. Soc. \textbf{346} (1994), no.~1, 1--62.

\bibitem{hochster:csopi}
\bysame, \emph{Comparison of symbolic and ordinary powers of ideals}, Invent.
  Math. \textbf{147} (2002), no.~2, 349--369.

\bibitem{huckaba:pihsad}
S.~Huckaba and C.~Huneke, \emph{Powers of ideals having small analytic
  deviation}, Amer. J. Math. \textbf{114} (1992), no.~2, 367--403.

\bibitem{huneke:tcaia}
C.~Huneke, \emph{Tight closure and its applications}, CBMS Conference 
  Ser.~in Math., vol.~88, Amer.~Math.~Soc., Providence, R.I., 1996.

\bibitem{huneke:ssrlr}
C.~Huneke and J.~Ribbe, \emph{Symbolic squares in regular local rings}, Math.
  Z. \textbf{229} (1998), 31--44.

\bibitem{kurano:pimsppi}
K.~Kurano and P.~C. Roberts, \emph{The positivity of intersection
  multiplicities and symbolic powers of prime ideals}, Compositio Math.
  \textbf{122} (2000), no.~2, 165--182.

\bibitem{matsumura:crt}
H.~Matsumura, \emph{Commutative ring theory}, second ed., Studies in
  Advanced Mathematics, vol.~8, University Press, Cambridge, 1989.

\bibitem{nagata:lr}
M.~Nagata, \emph{Local rings}, Interscience Tracts in Pure and Appl~Math.,
  No. 13, John Wiley \& Sons\, New York-London, 1962.

\bibitem{northcott:rilr}
D.~G. Northcott and D.~Rees, \emph{Reductions of ideals in local rings}, Proc.
  Cambridge Philos. Soc. \textbf{50} (1954), 145--158.

\bibitem{peskine:dpfcl}
C.~Peskine and L.~Szpiro, \emph{Dimension projective finie et cohomologie
  locale}, Inst. Hautes \'{E}tudes
  Sci. Publ. Math. (1973), no.~42, 47--119.

\bibitem{ratliff:qldafccpiII}
L.~J. {Ratliff, Jr.}, \emph{On quasi-unmixed local domains, the altitude
  formula, and the chain condition for prime ideals ({II})}, Amer. J. Math.
  \textbf{92} (1970), 99--144.

\bibitem{roberts:vimpc}
P.~C. Roberts, \emph{The vanishing of intersection multiplicities of perfect
  complexes}, Bull. Amer. Math. Soc. (N.S.) \textbf{13} (1985), no.~2,
  127--130.

\bibitem{roberts:rdsmcgpnc}
\bysame, \emph{Recent developments on {S}erre's multiplicity conjectures:
  {G}abber's proof of the nonnegativity conjecture}, Enseign. Math. (2)
  \textbf{44} (1998), 305--324.

\bibitem{sather:mdium}
S.~Sather-Wagstaff, \emph{Multiplicities and a dimension inequality for unmixed
  modules}, J. Algebra \textbf{238} (2001), no.~1, 372--388.

\bibitem{sather:dicmr}
\bysame, \emph{A dimension inequality for {C}ohen-{M}acaulay rings}, Trans.
  Amer. Math. Soc. \textbf{354} (2002), no.~3, 993--1005.

\bibitem{sather:isppi}
\bysame, \emph{Intersections of symbolic powers of prime ideals}, J. London
  Math. Soc. (2) \textbf{65} (2002), no.~3, 560--574.

\bibitem{serre:alm}
J.-P. Serre, \emph{Alg\`ebre locale - multiplicit\'es, troisi\`eme \'edition},
  Springer-Verlag, 1975.

\end{thebibliography}

\providecommand{\bysame}{\leavevmode\hbox to3em{\hrulefill}\thinspace}


\end{document}